\documentclass[12pt]{article}
\usepackage[cp1251]{inputenc} 
\usepackage[dvips]{graphicx}
\usepackage{amssymb}
\begin{document}
 \thispagestyle{empty}


\centerline{\textbf{ON THE GENERAL FORM OF LINEAR FUNCTIONAL ON}}
\centerline{\textbf{ THE  HARDY
SPACES $H^1$  OVER COMPACT ABELIAN GROUPS}}
\centerline{\textbf{ AND SOME OF ITS
APPLICATIONS}}
\vspace{1cm}
\centerline{\textbf{A. R. Mirotin}}

\section{Introduction}
\label{1}

 It is well known that the dual spaces $H^1$ and $BMO$ play a fundamental
role in PDEs and some other branches of analysis because they are often
the natural replacements for $L^1$ and $L^\infty$.

 Hardy spaces  $H^p$ over compact  Abelian groups with totally ordered duals were introduced by Helson and  Lowdenslager \cite{HL}, see also \cite[Chapter 8]{Rud}.
 In this note we generalize the celebrated Fefferman's  theorems on the dual of the complex and real Hardy  spaces  $H^1$ on the circle group (see \cite{F},  \cite{Gar})  to the case of arbitrary compact and connected Abelian group and give some applications of these results to lacunary 	multiple Fourier series,  multidimensional Hankel operators, and  atomic theory on the two dimensional torus. Our main tool in this study is the theory of Hilbert transform on (locally-)compact Abelian groups \cite{Rud}, \cite{ijpam}.

  In the following $G$ stands for compact  Abelian group  with the normalized Haar measure $m$ and totally ordered dual $X$, $X_+:=\{\chi\in X:\chi\geq {\bf 1}\}$  the positive cone in   $X$ (${\bf 1}$ denotes the unit character).
    As is well known, a (discrete) Abelian group $X$ can be totally ordered if and
only if it is torsion-free (see, for example, \cite{Rud}), which in turn is equivalent to the
condition that its character group $G$ is connected  \cite{Pont}; the total order on $X$
here is not, in general, unique.
 In applications, often $X$ is a dense
 subgroup of $\mathbb{R}^n$ endowed with the discrete topology so that $G$ is its Bohr compactification, or $X = \mathbb{Z}^n$ so that
$G = \mathbb{T}^n$ is the $n$-torus ($\mathbb{T}$ is the circle group and $\mathbb{Z}$ is the group of integers). For other examples we refer to \cite{SbMath}.

We denote by  $\widehat{\varphi}$
 the Fourier transform of $\varphi\in L^1(G)$, and by $\|\cdot\|_\infty$ the norm in $L^\infty(G)$. We put also
  $$
  \|f\|_p=\left(\int\limits_G|f|^pdm\right)^{1/p}
  $$
  for $f\in L^p(G)\ (0< p<\infty)$.

The\textit{ Hardy space} $H^p(G)\ (1\leq p\leq\infty)$ {\it over}
$G$ (with respect to the distinguished order on $X$)  is the subspace of  $L^p(G)$ defined as follows
$$
H^p(G)=\{f\in L^p(G):\widehat{f}(\chi)=0\ \forall\chi\notin X_+\}.
$$

 In particular $H^2(G)$ is the subspace of $L^2(G)$ with Hilbert basis $X_+$.
Let $P_+:L^2(G)\to H^2(G)$ be the orthogonal projection, $P_-=I-P_+$.

 For  every $u\in L^2(G,\mathbb{R})$
 there is a unique  $\widetilde{u}\in L^2(G,\mathbb{R})$ such that $\widehat{\widetilde{u}}({\bf 1})=0$ and
$u+{\it i}\widetilde{u}\in H^2(G)$.
The linear continuation of the mapping
$u\mapsto \widetilde{u}$ to the complex $L^2(G)$ is called a {\it Hilbert transform} on
$G$. This operator extends to a  bounded operator $\mathcal{H}:\varphi\mapsto \widetilde\varphi$ on $L^p(G)$ for $1<p<\infty$ (generalized Marcel Riesz's inequality), in particular  $\|\mathcal{H}\varphi\|_2\leq\|\varphi\|_2$ for every $\varphi\in L^2(G)$  \cite[8.7]{Rud}, \cite[Theorem 8, Corollary 20]{ijpam}.  Note also that  the Hilbert transform is a continuous map from $L^1(G)$ to $L^p(G)$ for $0<p<1$ (see, e.g.,\cite[Theorem 8.7.6]{Rud}).

\textbf{Definition 1} \cite{Trudy} (cf. \cite[p. 189]{Nik1}). We define the space $BMO(G)$ of
\textit{functions of bounded mean oscillation on   $G$} and its subspace $BMOA(G)$,
 as follows
$$
BMO(G):=\{f+\widetilde{g}: f,g\in L^{\infty}(G)\},
BMOA(G):=BMO(G)\cap H^1(G),
$$
$$
\|\varphi\|_{BMO}:=\inf\{\|f\|_{\infty}+\|g\|_{\infty}:\varphi=f+\widetilde{g},
f,g\in L^{\infty}(G)\} \ (\varphi\in BMO(G)).
$$

 \textbf{Lemma 1.}
\textit{The following equalities hold:}

$
(1)\  BMO(G)=P_-L^{\infty}(G)+P_+L^{\infty}(G),
$
\textit{with an equivalent norm}
$$
\|\varphi\|_{\ast}:=\inf\{\max(\|f_1\|_\infty, \|g_1\|_\infty):\varphi=P_-f_1+P_+g_1, f_1, g_1\in L^\infty(G)\};
$$

(2) \ $BMOA(G)=P_+L^{\infty}(G)$. \textit{Moreover, for the norm}
$$\|\varphi\|_{\ast}=\inf\{\|h\|_\infty:\varphi=P_+h,\
h\in L^\infty(G)\}
$$
\textit{in this space the following inequalities take place: \footnote{Here we correct a typo made in \cite[p. 139]{Trudy}.}}
$$\frac{2}{3}\|\varphi\|_
{BMO}\leq \|\varphi\|_
{\ast}\leq 2\|\varphi\|_
{BMO}.$$

\textbf{Proof}. The statement of the Lemma, except for the left hand side of the last inequality is contained in \cite[Proposition 3]{Trudy}. But the equality  $\varphi=P_+h,\
h\in L^\infty(G)$ implies $\varphi=f+\widetilde{g}$, where  $f=1/2(h+\widehat{h}({\bf 1})),\ g=i/2h$  \cite[Lemma 2]{Trudy}.
Thus
$\|\varphi\|_
{BMO}\leq\|f\|_\infty+\|g\|_\infty\leq 3/2\|h\|_\infty$,
and it left to go to the infimum when  $h$ runs through the space $L^\infty(G)$. $\Box$


\section{Duality theorems}
\label{2}

In the following we denote by $Y^*$ the dual of the Banach space $Y$.

{\bf Theorem 1.} \textit{
For every $\varphi\in BMOA(G)$ the formula
$$
F(f)=\int\limits_G f\overline \varphi dm\eqno(1)
$$
defines a linear functional on  $H^\infty(G)$,
and this functional extends uniquely to a continuous linear functional $F$  on $H^1(G)$. Moreover, the correspondence  $\varphi\mapsto F$ is an isometrical isomorphism of $(BMOA(G),\|\cdot\|_*)$  and $H^1(G)^*$, and a topological isomorphism of $(BMOA(G),\|\cdot\|_{BMO})$  and $H^1(G)^*$.}

 \textbf{Proof}. Let   $\varphi\in BMOA(G)$ and the functional  $F$ on  $H^\infty(G)$ is defined by the formula  (1). By Lemma 1,  $\varphi=P_+h$, where $h\in L^\infty(G)$. Moreover,
$$
F(f)=\int\limits_G f\overline{P_+h} dm=\int\limits_G P_+f\overline{h} dm=\int\limits_G f\overline{h} dm,
$$
 which implies that $|F(f)|\leq\|h\|_\infty\|f\|_1\ (f\in H^\infty(G))$. Passing to the infimum over $h$, we get $|F(f)|\leq\|\varphi\|_{\ast}\|f\|_1$ for every  $f\in H^\infty(G)$. Since $H^\infty(G)$ is dense in  $H^1(G)$ \cite[Lemma 1]{SbMath},  $F$ extends uniquely to a continuous linear functional $F$  on $H^1(G)$, and  $\|F\|\leq\|\varphi\|_{\ast}$.

 Conversely, for every linear functional $F\in H^1(G)^\ast$   there is a norm preserving extension $F$ to $L^1(G)$. Therefore there is such  $g\in L^\infty(G)$, that
$$
F(f)=\int\limits_G f\overline{g} dm\ (f\in L^1(G)),\ \mbox{ and }\ \|F\|=\|g\|_\infty.
$$
So for every  $f\in H^\infty(G)$ we have
$$
F(f)=\int\limits_G P_+f\overline{g} dm=\int\limits_G f\overline{P_+ g} dm.
$$
Thus   $F$ has the representation (1) with $\varphi=P_+ g\in BMOA(G)$ (Lemma 1). Moreover, $\|F\|=\|g\|_\infty\geq \|\varphi\|_{\ast}$. It follows   that $\|F\|= \|\varphi\|_{\ast}$.
We conclude that the linear map $(BMOA(G),\|\cdot\|_*) \to
H^1(G)^*,\ \varphi\mapsto F$ is surjective and isometric and as a result it is bijective.  Application of Lemma 1 completes the proof. $\Box$

 Note that by \textit{trigonometric polynomial on} $G$ we as usual mean the linear combination of characters (with complex coefficients in general).

\textbf{Definition 2}. We define the space $H^1_{\mathbb{R}}(G)$ (the \textit{real $H^1$ space on $G$})
as the completion of the space ${\rm Pol}(G,\mathbb{R})$ of real-valued trigonometric polynomials on $G$
 with respect to the norm
 $$
 \|q\|_{1\ast}:=\|P_-q\|_1+\|P_+q\|_1.
 $$
 We denote the norm in  $H^1_{\mathbb{R}}(G)$ by  $\|\cdot\|_{1\ast}$, too.

In the  next proposition we list  several impotent properties of $H^1_{\mathbb{R}}(G)$.

 \textbf{Proposition 1.} (i) \textit{Projectors $P_{\pm}$, and the Hilbert transform $\mathcal{H}$  are  bounded operators on } $H^1_{\mathbb{R}}(G)$;

 (ii) \textit{restrictions $P_{\pm}|{\rm Pol}(G,\mathbb{R})$ extend to  bounded operators $P_{\pm}^1$ from $ H^1_{\mathbb{R}}(G)$ to $L^1(G)$ and
 $$
 \|f\|_{1*}=\|P_-f\|_{1*}+\|P_+f\|_{1*}=\|P_-^1f\|_{1}+\|P_+^1f\|_{1}\    (f\in H^1_{\mathbb{R}}(G));
 $$
 }

 (iii) $ H^1_{\mathbb{R}}(G)={\rm Im}P_-\dotplus {\rm Im}P_+$ (\textit{the direct sum of closed subspaces});

 (iv) $\cup_{p>1}L^p(G,\mathbb{R})\subset  H^1_{\mathbb{R}}(G) \subset L^1(G,\mathbb{R})$;

 (v) $\|f\|_{\mathcal{H}}:=\|f\|_1+\|\mathcal{H}f\|_1$ \textit{is an equivalent norm in} $ H^1_{\mathbb{R}}(G)$;

 (vi)  $H^1_{\mathbb{R}}(G)={\rm Re}H^1(G)$.

\textbf{Proof}. (i) The boundedness of $P_{\pm}$ follows from the inequalities  $\|P_{\pm} q\|_{1\ast}=\|P_{\pm} q\|_{1}\leq \|q\|_{1\ast}$, and
the boundedness of the Hilbert transform is the consequence of the equality $i\mathcal{H}q=2P_+q-q-2\widehat q({\bf 1})$ \cite[Lemma 22]{ijpam}, since $\|\widehat q({\bf 1})\|_{1*}=|\widehat q({\bf 1})|\leq\|q\|_1\leq \|q\|_{1*}\ (q\in {\rm Pol}(G,\mathbb{R}))$.

(ii) Inequalities $\|P_{\pm}q\|_1\leq \|q\|_{1\ast}$ implies that $P_{\pm}$ extend to  bounded operators $P_{\pm}^1$ from $ H^1_{\mathbb{R}}(G)$ to $L^1(G)$. The first equality  follows from (i) and the equality $ \|q\|_{1\ast}=\|P_-q\|_{1\ast}+\|P_+q\|_{1\ast}\ (q\in {\rm Pol}(G,\mathbb{R}))$.

Now we claim that there is a continuous embedding $H^1_{\mathbb{R}}(G)\subset L^1(G,{\mathbb{R}})$. Indeed, the norms $\|\cdot\|_1$ and $\|\cdot\|_{1*}$ in ${\rm Pol}(G,\mathbb{R})$ are comparable, since $\|q\|_1\leq \|q\|_{1*}$ for $q\in {\rm Pol}(G,\mathbb{R})$. Moreover, they are compatible in the sense that  every sequence $(q_n)$, which is fundamental with respect to both norms and converges to the zero element  with respect to $\|\cdot\|_1$, also converges to the zero element  with respect to $\|\cdot\|_{1*}$ \cite[p. 13]{GS}. For the proof of this statement first note that $\widehat {q_n}\to 0$ (the uniform convergence on $X$). Since
$$
\|P_-(q_n-q_m)\|_1+\|P_+(q_n-q_m)\|_1= \|q_n-q_m\|_{1*}\to 0\ (m,n\to\infty),
 $$
 $\|P_{\pm}q_n- h_{\pm}\|_1\to 0$ for some $h_{\pm}\in L^1(G, \mathbb{R})$, and therefore $\widehat {P_{\pm}q_n}\to \widehat {h_{\pm}}$ (the uniform convergence on $X$). On the other hand $2P_+q_n=i\mathcal{H}q_n+q_n+2\widehat q_n({\bf 1})$ and therefore
 $$
2\widehat {P_+q_n}=i\widehat {\mathcal{H}q_n}+\widehat {q_n}+2\widehat q_n({\bf 1})=
(i{\rm sgn}_{X_+}+1)\widehat {q_n}+2\widehat q_n({\bf 1})\to 0 (n\to \infty)
$$
(see \cite[Lemma 5]{ijpam}). So $h_{\pm}=0, \|q_n\|_{1*}\to 0$ and the  continuous embedding $H^1_{\mathbb{R}}(G)\subset L^1(G,{\mathbb{R}})$ follows \cite[p. 14]{GS}.

 Hence
for $f\in H^1_{\mathbb{R}}(G)$ and a sequence $q_n\in {\rm Pol}(G,\mathbb{R})$  such that $\|q_n-f\|_{1*}\to 0$ we have  $\|q_n-f\|_{1}\to 0$ and thus
$$
\|f\|_{1*}=\lim_n\|q_n\|_{1*}=\lim_n(\|P_-q_n\|_1+\|P_+q_n\|_1)=\|P_-^1f\|_{1}+\|P_+^1f\|_{1}.
$$

(iii) Since $P_-, P_+$ are bounded projectors on  $H^1_{\mathbb{R}}(G)$, their images ${\rm Im}P_-$, and ${\rm Im}P_+$ are closed subspaces of $H^1_{\mathbb{R}}(G)$. The equality
$ H^1_{\mathbb{R}}(G)={\rm Im}P_-+ {\rm Im}P_+$ follows from the boundedness of $P_{\pm}$. Now if $f\in {\rm Im}P_-\cap {\rm Im}P_+$, the property (ii) implies that $\|f\|_{1*}=0$.

(iv) The right hand inclusion was proved above (see the proof of (ii)). Let $f\in L^p(G,\mathbb{R})\ (p>1)$ and $\|f-q_n\|_p\to 0, n\to\infty (q_n\in {\rm Pol}(G,\mathbb{R}))$. It is known  that the map $P_+$ can be extended from ${\rm Pol}(G,\mathbb{R})$ to bounded projection $\Phi:L^p(G)\to H^p(G)$  \cite[8.7.2]{Rud}. Then
$$
\|P_+q_n-\Phi f\|_1+\|P_-q_n-(I-\Phi) f\|_1=\|\Phi(q_n- f)\|_1+\|(I-\Phi)(q_n- f)\|_1\leq
$$
$$
\|\Phi(q_n- f)\|_p+\|(I-\Phi)(q_n- f)\|_p\to 0\ (n\to\infty).
$$
Thus $P_+q_n\to \Phi f$ and $P_-q_n\to (I-\Phi) f$ in $L^1(G)$. This implies that
$$
\|q_n-q_m\|_{1*}:=\|P_-(q_n-q_m)\|_1+\|P_+(q_n-q_m)\|_1\to 0\ (n,m\to\infty),
$$
and therefore $\|q_n-g\|_{1*}\to 0$ for some $g\in H^1_{\mathbb{R}}(G)$. But from (ii)
it follows that  $P^1_{\pm}q_n\to P^1_{\pm}g$ in $L^1$ norm. So $P^1_{+}g=\Phi f, P^1_{-}g=(I-\Phi)f$,
and we conclude that $f=P^1_{-}g+P^1_{+}g=g\in H^1_{\mathbb{R}}(G)$.

(v) Equalities $P_{\pm}q=1/2(\pm i\mathcal{H}q +q)\pm \widehat q({\bf 1})$ show that
$\|P_{\pm}q\|_1\leq 3/2(\|q\|_1+\|\mathcal{H}q\|_1$ and therefore $\|q\|_{1*}\leq 3\|q\|_{\mathcal{H}}$. On the other hand $\|q\|_{\mathcal{H}}\leq \|q\|_{1*}+\|\mathcal{H}q\|_{1*}\leq (1+\|\mathcal{H}\|)\|q\|_{1*}$,
since $\mathcal{H}$ is bounded on  $H^1_{\mathbb{R}}(G)$.

(vi) Let $u\in H^1_{\mathbb{R}}(G)$ and $\|q_n-u\|_{1*}\to 0$ for some $q_n\in {\rm Pol}(G,\mathbb{R})\ (n\to\infty)$. In view of (v) $\|q_n-u\|_{1}\to 0$ and $\|\mathcal{H}q_n-\mathcal{H}u\|_{1}\to 0 \ (n\to\infty)$. Since,  by the definition of the Hilbert transform,  $q_n+i\mathcal{H}q_n\in H^2(G)\subset H^1(G)$, it follows that $u+i\mathcal{H}u\in H^1(G)$ and thus $u\in {\rm Re}H^1(G)$.

 Conversely, let $u\in {\rm Re}H^1(G)$ and $f=u+iv\in H^1(G)$. Then $\|p_n-f\|_1\to 0 \ (n\to\infty)$ for some $p_n\in {\rm Pol}(G)$ \cite[Lemma 1]{SbMath}. Adding, if necessary, to $f$ a pure imaginary constant, we can assume that $\widehat{v}(\textbf{1})=0$.  Let $p_n=q_n+ih_n$, where $q_n, h_n\in {\rm Pol}(G,\mathbb{R})$. Then  $\|q_n-u\|_1\to 0, \|h_n-v\|_1\to 0 $, and $\widehat{h_n}(\textbf{1})\to \widehat{v}(\textbf{1})=0\ (n\to\infty)$. Replacing, if necessary, $p_n$ with $p_n-i\widehat{h_n}(\textbf{1})$ we can assume that $\widehat{h_n}(\textbf{1})=0$.  Fix $p\in(0,1)$. Since  the Hilbert transform continuously maps $L^1(G)$ into $L^p(G)$ \cite[Chapter 8]{Rud},   we have $\|\mathcal{H}q_n-\mathcal{H}u\|_p\to 0$. By the definition of the Hilbert transform  $\mathcal{H}q_n=h_n$,  which implies that $\|h_n-\mathcal{H}u\|_p\to 0$. On the other hand, $\|h_n-v\|_p\leq\|h_n-v\|_1 \to 0$ and therefore $v=\mathcal{H}u$. It follows that
 $$
 \|q_n-u\|_{\mathcal{H}}=\|q_n-u\|_1+\|\mathcal{H}q_n-\mathcal{H}u\|_1=\|q_n-u\|_1+ \|h_n-v\|_1\to 0 \ (n\to\infty),
 $$
and the application of the statement (v) finishes the proof. $\Box$

 Now we are in position to prove the real version of  Feffermans' duality theorem. By $BMO(G,\mathbb{R})$ we denote the subspace of real-valued functions from $BMO(G)$.

\textbf{ Theorem 2.} \textit{
For every  $\varphi\in BMO(G,\mathbb{R})$ the linear functional
$$
F(q)=\int\limits_G q \varphi dm \eqno(2)
$$
 on  ${\rm Pol}(G,\mathbb{R})$ extends uniquely to a continuous linear functional $F$  on $H^1_{\mathbb{R}}(G)$. Moreover, the correspondence  $\varphi\mapsto F$ is an isometrical isomorphism of $(BMO(G,\mathbb{R}),\|\cdot\|_*)$  and $H^1_{\mathbb{R}}(G)^*$, and  a topological isomorphism of  $(BMO(G,\mathbb{R}),\|\cdot\|_{BMO})$  and $H^1_{\mathbb{R}}(G)^*$.}

\textbf{Proof}.
 Let  $\varphi\in BMO(G,\mathbb{R})$ and the functional  $F$ on  ${\rm Pol}(G,\mathbb{R})$ is defined by the formula  (2). By Lemma 1,  $\varphi=P_-g+P_+h$, where $g, h\in L^\infty(G)$.
 Then for every $q\in {\rm Pol}(G,\mathbb{R})$ we have
$$
F(q)=\int\limits_G P_-g\overline{q} dm+\int\limits_G P_+h\overline{q} dm=\int\limits_G g\overline{P_-q} dm+\int\limits_G h\overline{P_+q} dm,
$$
 which implies that $|F(q)|\leq\max(\|g\|_\infty,\|h\|_\infty)(\|P_-q\|_1+\|P_+q\|_1)$.
So, $\|F\|\leq \|\varphi\|_*$, and there is
a unique extension of $F$ to a continuous linear functional $F$  on $H^1_{\mathbb{R}}(G)$ with the same norm.

 Conversely, for every linear functional $F\in H^1_{\mathbb{R}}(G)^\ast$ let $F_-$ and $F_+$ denote its restrictions to ${\rm Im}P_-$ and ${\rm Im}P_+$  respectively. If $f\in {\rm Im}P_{\pm}, f=P_{\pm}g$, where $g\in H^1_{\mathbb{R}}(G)$, then $|F_{\pm}(f)|\leq \|F\|\|P_{\pm}g\|_{1*}=\|F\|\|f\|_1$.     Therefore these functionals  extend to linear bounded functionals $F_-$ and $F_+$ on  $L^1(G,\mathbb{R})$ with preservation of norms. Let $g_{\pm}\in L^\infty(G)$ be such that
$$
F_{\pm}(f)=\int\limits_G f \overline{g_{\pm}} dm\ (f\in L^1(G,\mathbb{R})), \mbox{ and } \|F_{\pm}\|=\|g_{\pm}\|_\infty.
$$
It follows that
$$
F(q_{\pm})=\int_G P_{\pm}q_{\pm} \overline{g_{\pm}} dm=\int_G q_{\pm} \overline{P_{\pm}g_{\pm}} dm
$$
 for $q_{\pm}\in P_{\pm}{\rm Pol}(G,\mathbb{R})$.
Since every $q\in  {\rm Pol}(G,\mathbb{R})$ has the form $q=q_+ + q_-$ where $q_{\pm}\in P_{\pm}{\rm Pol}(G,\mathbb{R})$ we have
$$
F(q)=F(q_+)+F(q_-)=\int\limits_G q_{+} \overline{P_{+}g_{+}} dm+\int\limits_G q_{-} \overline{P_{-}g_{-}} dm=\int\limits_G q (\overline{P_{+}g_{+}} +\overline{P_{-}g_{-}})dm.
$$
Thus if we put $\varphi:=\overline{P_{+}g_{+}} +\overline{P_{-}g_{-}}$, the equality (2) holds.
Since $F(q)$ is real-valued for every $q\in  {\rm Pol}(G,\mathbb{R})$ so is $\varphi$. In fact,  putting   $q=\chi+\overline{\chi}$ and then   $q=i(\chi-\overline{\chi})$,
we deduce from the equality $F(q)=\overline{F(q)}$ in view of  (2) that for every $\chi\in X$
$$
\int\limits_G(\varphi-\overline\varphi)\chi dm\pm\int\limits_G(\varphi-\overline\varphi)\overline\chi dm=0.
$$
The last two equalities imply  that the Fourier transform of  the  imaginary part of $\varphi$ equals to zero. Thus $\varphi=P_{+}g_{+} +P_{-}g_{-}\in BMO(G,\mathbb{R})$. Moreover, $\|F\|\geq \|F_{\pm}\|=\|g_{\pm}\|_\infty$.  Passing to the infimum over $g_{\pm}$, we get $\|F\|\geq \|\varphi\|_*$, and the rest of the proof is exactly the same as the rest of the proof of  Theorem 1.  $\Box$

\section{Applications. Lacunary series}
\label{3}

Now we apply  Theorem 1 in order to generalize  some results on lacunary series  in one variable (see \cite[p. 191, 1.6.4  (a), (d)]{Nik1}).
The following definition is implicitly contained in \cite[8.6]{Rud}.

\textbf{Definition 3}. We call a subset $E\subset X_+$   \textit{lacunary} (in the sense of Rudin) if there is a constant  $K=K_E$ such that  the number of terms of  the set $\{\xi\in E:\chi\leq\xi\leq\chi^2\}$ do not exceed $K$ for every $\chi\in X_+$.

\textbf{Theorem 3.} \textit{Let  $\varphi\in L^1(G)$, and $\widehat \varphi$ vanishes outside some lacunary set $E$. Then  $\varphi\in BMOA(G)$ if and only if  $\varphi\in H^2(G)$; moreover,}
$$
\|\varphi\|_2\leq\|\varphi\|_{BMO}\leq 3\sqrt{K_E}\|\varphi\|_2.
$$

\textbf{Proof}. By the Lemma 1 $BMOA(G)\subset H^2(G)$.
Now let  $\varphi\in H^2(G)$ and $\widehat \varphi$ vanishes outside $E\subset X_+$. Then $\widehat \varphi\in l_2(E)$. We claim that the functional
$$
\Lambda(f):=\sum_{\chi\in E}\widehat{f}(\chi)\overline{\widehat{\varphi}(\chi)}
$$
is defined and bounded  on  $H^1(G)$. Indeed, for every  $f\in H^1(G)$ the restriction of $\widehat{f}$ to $E$ belongs to  $l_2(E)$ by \cite[Theorem 8.6]{Rud}. Moreover,  $\|\widehat{f}\|_{l_2(E)}\leq 2\sqrt{K_E}\|f\|_1$ by \cite[p. 214, (5)]{Rud} and therefore
$$
|\Lambda(f)|\leq \|\widehat{f}\|_{l_2(E)}\|\widehat{\varphi}\|_{l_2(E)}\leq 2\sqrt{K_E}\|\widehat{\varphi}\|_{l_2(E)}\|f\|_1.
$$
Hence accordingly to Theorem 1 there exists such  $\varphi_1\in BMOA(G)$ that
$$
\Lambda(f) = \int\limits_G f\overline{\varphi_1} dm \ (f\in H^\infty(G)).\eqno(3)
$$

On the other hand it follows from the Plancherel Theorem  that
$$
\Lambda(f)=\sum_{\chi\in X}\widehat{f}(\chi)\overline{\widehat{\varphi}(\chi)}=\int\limits_G f\overline \varphi dm\eqno(4)
$$
 for every  $f\in H^\infty(G)$.
 Now for $f\in X_+$ formulas  (3) and (4) imply that the function $\widehat{\varphi- \varphi_1}$ vanishes on    $X_+$. Since this function is concentrated on   $X_+$, we get  $\varphi = \varphi_1$, which completes the proof of the first statement.

 Next, let $\varphi=f+\widetilde{g}$ where $f,g\in L^{\infty}(G)$. Since $\|\widetilde{g}\|_2\leq\|g\|_2$ \cite[Theorem 8]{ijpam}, we have
$$
\|\varphi\|_2\leq\|f\|_2+\|\widetilde{g}\|_2\leq \|f\|_\infty+\|g\|_\infty,
$$
and the first inequality follows. On the other hand
if the functional $F$ is defined by the formula  (1), then $\|F\|=\|\varphi\|_{\ast}$ (see the proof of Theorem  1). But (4) implies that  $F=\Lambda$, and we already have seen that
$$
\|\Lambda\|\leq 2\sqrt{K_E}\|\widehat{\varphi}\|_{l_2(E)}=2\sqrt{K_E}\|\varphi\|_2.
  $$
  The application of Lemma 1 completes the proof. $\Box$

 For the next corollary of Theorem 1 recall  that \textit{Hankel operator}  $H_{\varphi}:H^2(G)\rightarrow  H^2_-(G)$ ($ H^2_-(G)=L^2(G)\ominus H^2(G)$) with symbol $\varphi\in L^2(G)$ is initially defined on the subspace of trigonometric polynomials of analytic type (the linear span of $X_+$) by the equality
 $$
 H_{\varphi}f=P_-({\varphi}f),
 $$
  (see, e.g., \cite{Trudy}).

 For compact Abelian groups the fundamental  Nehari Theorem for Hankel forms  was proved by J.~ Wang \cite{Wang}. Its version for  Hankel operators looks as follows.

\textbf{Theorem} (Z. Nehari). \textit{A bounded operator $H:H^2(G)\to  H^2_-(G)$ is of the form $H_\varphi$ for some $\varphi\in L^\infty(G)$  if and only if
$$
HS_\chi=P_-S_\chi H\  \forall \chi\in X_+,\eqno(5)
$$
where $S_\chi f:=\chi f\ (f\in L^2(G))$. Moreover, $\|H\|=\|\varphi\|_\infty$ for some $\varphi\in L^\infty(G)$ such that $H=H_\varphi$.}

For the proof  one can  consider the bilinear form $(f,g)\mapsto\langle Hf,\overline g\rangle$ on $H^2(G)\times H^2(G)$ and apply the result from \cite{Wang} to  the corresponding bilinear Hankel form on $l_2(X_+)\times l_2(X_+)$ ($l_2(X_+)$ is isomorphic to $H^2(G)$ via the Fourier transform; see \cite{Dyba} for details). $\Box$

\textbf{Theorem 4.} \textit{Let $E$ be a lacunary set and ${\bf 1}\notin E$. Then for every  $\varphi\in H^2(G)$ such that  $\widehat \varphi$ vanishes outside   $E$ the operator  $H_{\overline{\varphi}}$ is bounded and }
$$
A_E\|\varphi\|_2\leq\|H_{\overline{\varphi}}\|\leq 6\sqrt{K_E}\|\varphi\|_2\eqno(6)
$$
\textit{where $A_E$  is independent of $\varphi$.}

\textbf{Proof}. By virtue of Theorem 3 and Lemma 1 $\varphi=P_+h, h\in L^\infty (G)$. So $P_-(\bar\varphi-\bar h)=0$, and therefore $H_{\bar\varphi}=H_{\bar h}$. Consequently $H_{\overline{\varphi}}$ is bounded and $\|H_{\overline{\varphi}}\|\leq \|h\|_\infty$. It follows that $\|H_{\overline{\varphi}}\|\leq \|\varphi\|_{\ast}$, and  Lemma 1 and Theorem 3 entail that
$\|\varphi\|_{\ast}\leq 6\sqrt{K_E}\|\varphi\|_2$. This proves the second inequality (without assumption that ${\bf 1}\notin E$).

To prove the first one consider  the following subspace of $H^2(G)$:
$$
H^2_E(G):=\{\varphi\in H^2(G): \widehat \varphi \mbox{ vanishes outside } E\}.
$$
Being isomorphic to $l_2(E)$ via the Fourier transform, this space is complete with respect to the $L^2$ norms. It is easy to verify that the seminorm $\|\varphi\|_H:=\|H_{\overline{\varphi}}\|$ is a norm in $H^2_E(G)$ if ${\bf 1}\notin E$. We claim that $H^2_E(G)$ is complete with respect to this norm as well. Indeed, let the sequence $(\varphi_n)\subset H^2_E(G)$ be fundamental with respect to  $\|\cdot\|_H$. Then $\|H_{\overline{\varphi_n}}-H\|\to 0 \ (n\to\infty)$ for some bounded operator $H, H:H^2(G)\to  H^2_-(G)$. But it is easy to verify that
$H_{\overline{\varphi_n}}S_\chi\xi=P_-S_\chi H_{\overline{\varphi_n}}\xi\  \forall \xi,\chi\in X_+$. Since the operator $H_{\overline{\varphi_n}}$ is bounded,  it follows that it satisfies the condition (5) and therefore $H$ satisfies (5), too. So by  Nehari's Theorem $H=H_g$ for some $g\in L^\infty(G)$. Let $f:=\overline{P_-g}$. Then $\widehat{\overline f}$ vanishes outside $X_-$
and $\|H_{\overline{\varphi_n}-\overline f}\|\to 0 \ (n\to\infty)$. Again by  Nehari's Theorem there is such $\varepsilon_n\in L^\infty(G)$ that $H_{\overline{\varphi_n}-\overline f}=H_{\overline{\varepsilon_n}}$ and $\|H_{\overline{\varphi_n}-\overline f}\|=\|\varepsilon_n\|_\infty$. Put $\psi_n:=\overline{\varphi_n}-\overline f-\overline{\varepsilon_n}$. Then $\psi_n\in H^2(G)$, since $H_{\psi_n}=0$. Thus for every $\chi\in X_+\setminus E, \chi\ne \textbf{1}$ we have
$$
\widehat f(\chi)=\widehat{\varphi_n}(\chi)-\widehat{\overline{\psi_n}}(\chi)-\widehat{\varepsilon_n}(\chi)=-\widehat{\varepsilon_n}(\chi),
$$
and therefore $|\widehat f(\chi)|\leq \|\varepsilon_n\|_\infty$. It follows that $f\in H^2_E(G)$ and since $\|\varphi_n-f\|_H\to 0\  (n\to\infty)$  the space $(H^2_E(G), \|\cdot\|_H)$ is complete.
 In view of the second inequality  in (6) and the well known Banach Theorem, the norms $\|\cdot\|_2$ and $\|\cdot\|_H$ are equivalent in $H^2_E(G)$. $\Box$

\textbf{Remark 1.} If we take $G=\mathbb{T}^n$, the $n$-dimensional torus, and choose some linear order on its dual group $\mathbb{Z}^n$ Theorems 3 and 4  turn into results on lacunary 	multiple Fourier series and  multidimensional Hankel operators. A description of
all linear orders on $\mathbb{Z}^n$  one can find in \cite{Teh}, \cite{Zajt}. The case of infinite dimensional torus  $\mathbb{T}^\infty$ (see, e.g., \cite[Examples 2, 3]{SbMath}) is also of interest.

\section{Some results related to  atomic theory on  $H^1_{\mathbb{R}}(\mathbb{T}^2)$}
\label{4}

The "atomic" theory for functions from $H^1_{\mathbb{R}}(\mathbb{T})$ was developed in  \cite{Coif}. A general approach to atomic decompositions was proposed by Coifman and Weiss  \cite{CW}, but their notion of an atom  \cite[p. 591]{CW} differs from ours (see  Definition 4 below; a remarkable feature of our atoms  is that these atoms have only partial cancellation conditions), and the definition of Hardy  spaces  in \cite[p. 592]{CW} differs from ours, too. For more resent results in this area see, e.g., \cite{Bow} --- \cite{Dek}.

The problem of developing an atomic theory for Hardy spaces on the polydisc was posed in \cite[p. 642]{CW}.
  In this section, we  get some results related to an atomic theory for  $H^1_{\mathbb{R}}(\mathbb{T}^2)$.  It should be noted that we consider the
   last  space with respect to the lexicographic order on the dual group $\mathbb{Z}^2$  of $\mathbb{T}^2$.

\textbf{Definition 4}.  By a  $\mathbb{T}^2$-\textit{atom} we  mean  either the function $\textbf{1}$ or a real-valued function $a(\theta_1,\theta_2)$ supported on a rectangle
$J_1\times J_2\subseteq \mathbb{T}^2$ having the property

(i) $|a(\theta_1,\theta_2)| \leq \min\{1/|J_1|, 1/|J_2|\}$,

${\rm (ii)} \int\limits_{J_1} a(\theta_1,\theta_2)d\theta_1 = 0 =\int\limits_{J_2} a(\theta_1,\theta_2)d\theta_2$
 for every $(\theta_1,\theta_2)\in J_1\times J_2$.

By  a  $\mathbb{T}^1$-\textit{atom} we  mean  either the function $\textbf{1}$ or a real-valued function $a(\theta_1)\ (a(\theta_2))$ on $\mathbb{T}^2$ supported on a rectangle
$J_1\times  \mathbb{T}$  (respectively $\mathbb{T}\times J_2$)  having the property

(${\rm i}^\prime$)  $|a(\theta_i)| \leq 1/|J_i|$,

(${\rm ii}^\prime$)  $\int\limits_{J_i} a(\theta_i)d\theta_i = 0$.

\noindent
Above $J_i$ denotes an arc in  $\mathbb{T}$ with normalized Lebesgue  measure $|J_i| \ (|\mathbb{T}|=1),\ i=1,2$.

\textbf{Proposition 2.}  \textit{Atoms form a bounded subset of $H^1_{\mathbb{R}}(\mathbb{T}^2)$ (and   generate a dense subspace of this space).}

We need two  preliminary results to prove Proposition 2.

 \textbf{Proposition 3.} \textit{The Hilbert transform in $L^2(\mathbb{T}^2)$ (with respect to the lexicographic order on  $\mathbb{Z}^2$) has the form }
$$
 \mathcal{H}f(t_1,t_2)
  =P.V.\int\limits_{-\pi}^{\pi}f(\theta_1,t_2)\cot\frac{t_1-\theta_1}{2}\frac{d\theta_1}{2\pi}
  $$
  $$
\hspace{18mm}  + \int\limits_{-\pi}^{\pi} P.V.\int\limits_{-\pi}^{\pi}f(\theta_1,\theta_2)\cot\frac{t_2-\theta_2}{2}\frac{d\theta_2}{2\pi}\frac{d\theta_1}{2\pi}
 $$
 $$
 \hspace{22mm}=\mathcal{H}_{\theta_1\to t_1}f(\theta_1,t_2) +\int\limits_{-\pi}^{\pi}(\mathcal{H}_{\theta_2\to t_2}f(\theta_1,\theta_2))\frac{d\theta_1}{2\pi} \eqno(7)
 $$
\textit{where $\mathcal{H}_{1}$ and $\mathcal{H}_{2}$ stand for the Hilbert transform in $L^2(\mathbb{T})$ in the first  and  the second variable independently}.

 \textbf{Proof}. First note that both summands in the right-hand side of formula (7) are continuous linear operators on   $L^2(\mathbb{T}^2)$. And as was mentioned in the introduction of this paper the  left-hand side of  (7) is continuous, too. So it remains to verify (7) for  functions of the form $f=u_i\otimes u_j$ where $u_i$ and $u_j$ run over some orthogonal base of  $L^2(\mathbb{T})$.
 To this end note that the equality $g=\mathcal{H}f$ is equivalent to
 $$
 \widehat g=-i{\rm sgn}_{X_+}\widehat f,
 $$
 where $X_+$ is the positive cone in $\mathbb{Z}^2$ with respect to the lexicographic order, and
 ${\rm sgn}_{X_+}(n_1,n_2)  := 1 (-1)$ for $(n_1,n_2)\in X_+\setminus\{0\}$ (respectively $(n_1,n_2)\notin X_+$), ${\rm sgn}_{X_+}(0,0)  :=0$   \cite{ijpam}. It is easy to verify that
 $$
 {\rm sgn}_{X_+}(n_1,n_2)  ={\rm sgn}_{\mathbb{Z}_+}(n_1)+{\rm sgn}_{\{0\}\times\mathbb{Z}_+}(n_1,n_2).\eqno(8)
 $$
 Since $\widehat f(n_1,n_2)=\widehat{u_i}(n_1)\widehat{u_j}(n_2)$  for all $(n_1,n_2)\in \mathbb{Z}^2$, we have using (8)
 $$
  \widehat g(n_1,n_2)=-i{\rm sgn}_{X_+}(n_1,n_2)\widehat f (n_1,n_2)
  $$
  $$
 = (-i{\rm sgn}_{\mathbb{Z}_+}(n_1)\widehat{u_i}(n_1))\widehat{u_j}(n_2)+ \widehat {u_i}(0)(-i{\rm sgn}_{\mathbb{Z}_+}(n_2)\widehat{u_j}(n_2)).
 $$
 This implies that
 $$
 g(t_1,t_2)=(\mathcal{H}u_i(t_1))u_j(t_2)+\widehat {u_i}(0)\mathcal{H}u_j(t_2)
 $$
 $$
= P.V.\int\limits_{-\pi}^{\pi}f(\theta_1,t_2)\cot\frac{t_1-\theta_1}{2}\frac{d\theta_1}{2\pi}+
 \int\limits_{-\pi}^{\pi} P.V.\int\limits_{-\pi}^{\pi}f(\theta_1,\theta_2)\cot\frac{t_2-\theta_2}{2}\frac{d\theta_2}{2\pi}\frac{d\theta_1}{2\pi},
 $$
concluding the proposition. $\Box$

 \textbf{Lemma 2.} \textit{There is  a universal constant $C>0$ such that $\|a\|_{1*}\leq C$ for every
  $\mathbb{T}^2$- or  $\mathbb{T}^1$-atom $a$.}

  \textbf{Proof}. We shall use the statement (v) of Proposition 1. It is evident that  $\|a\|_{1}\leq 1$ for every
  $\mathbb{T}^2$- or  $\mathbb{T}^1$-atom $a$.

  By Proposition 3 $\mathcal{H}a(t_1,t_2)=A_1(t_1,t_2)+A_2(t_2)$, where
  $$
  A_1(t_1,t_2) = \mathcal{H}_{1}a(t_1,t_2)= P.V.\int\limits_{-\pi}^{\pi}a(\theta_1,t_2)\cot\frac{t_1-\theta_1}{2}\frac{d\theta_1}{2\pi},
  $$
  and
  $$ A_2(t_2) =
 \int\limits_{-\pi}^{\pi} \mathcal{H}_{2}a(\theta_1,t_2) \frac{d\theta_1}{2\pi}.
  $$

 Let $a$ be a $\mathbb{T}^1$-atom, $a\ne\textbf{ 1}$. If $a$ depends of $\theta_1$ only, then $\|A_1\|_1=\|\mathcal{H}_{1}a\|_{L^1(dt_1/2\pi)}\leq{\rm const}$
by the classical result on  $\mathbb{T}^1$-atoms (see, e. g. \cite[p. 27]{Hoep}), and if $a$ depends of $\theta_2$ only, we get $A_1=0$.

Now let $a$ be a $\mathbb{T}^2$-atom, $a\ne\textbf{ 1}$. If $J_1$ is the arc appearing in (i) and
(ii), we can assume, without loss
of generality, that  $J_1 = (-\delta; \delta)$. As in the one dimensional case we can assume also, that $\delta<1/2$, since $|a|\leq 1$ for $\delta\geq 1/2$. Then we have
$$
\|A_1\|_1=\int\limits_{J_2}\left(\int\limits_{-2\delta}^{2\delta}|A_1(t_1,t_2)|\frac{dt_1}{2\pi}+ \int\limits_{2\delta<|t_1|<\pi}|A_1(t_1,t_2)|\frac{dt_1}{2\pi}\right)\frac{dt_2}{2\pi}
$$
$$
=\int\limits_{J_2}\left(I_1(t_2)+I_2(t_2)\right)\frac{dt_2}{2\pi}.
$$
Using the H$\ddot {\rm o}$lder's and  generalized Marcel Riesz's inequality, and property (i) of atom we obtain (below we denote by $C$ or ${\rm const}$ any universal constant)
$$
I_1(t_2)=\int\limits_{-2\delta}^{2\delta}|A_1(t_1,t_2)|\frac{dt_1}{2\pi}\leq C\sqrt{\delta}\|A_1(\cdot,t_2)\|_{L^2(\frac{dt_1}{2\pi})}
$$
$$
=C\sqrt{\delta}\|\mathcal{H}_{1}a(t_1, t_2)\|_{L^2(\frac{dt_1}{2\pi})}
\leq C\sqrt{\delta}\|a(\cdot, t_2)\|_{L^2(\frac{d\theta_1}{2\pi})}\leq \frac{C}{|J_2|}.
$$
Therefore, $\int_{J_2}I_1(t_2)dt_2/2\pi\leq {\rm const}.$

To estimate $I_2$ we use  the following  classical estimate for the one-dimensional case (see, e. g., \cite[p. 28]{Hoep})
$$
|A_1(t_1,t_2)|=|\mathcal{H}_{1}a(t_1, t_2)|\leq C\delta \|a(\cdot, t_2)\|_{L^1(\frac{d\theta_1}{2\pi})}t_1^{-2}\  (|t_1|\geq 2\delta). \eqno(9)
$$
It follows that
$$
I_2(t_2)=\int\limits_{2\delta<|t_1|<\pi}|A_1(t_1,t_2)|\frac{dt_1}{2\pi}\leq C\delta\int\limits_{J_1}|a(\theta_1, t_2)|\frac{d\theta_1}{2\pi}\int\limits_{2\delta<|t_1|<\pi}t_1^{-2}\frac{dt_1}{2\pi}\leq {\rm const},
$$
and therefore $\int_{J_2}I_2(t_2)dt_2/2\pi\leq {\rm const}$, as well. We conclude that
$$
\|A_1\|_1=\int\limits_{-\pi}^{\pi}\int\limits_{-\pi}^{\pi}|\mathcal{H}_{1}a(t_1, t_2)|\frac{dt_1}{2\pi}\frac{dt_2}{2\pi}\leq {\rm const}.\eqno(10)
$$
As regards $A_2$, the case when  $a$ is a $\mathbb{T}^1$-atom can be considered as above.
Let $a$ be a $\mathbb{T}^2$-atom, $a\ne \textbf{1}$. Then by the Fubini's theorem we get in view of (10) that
$$
\|A_2\|_1=\|A_2\|_{L^1(\frac{dt_2}{2\pi})}\leq\int\limits_{-\pi}^{\pi}\int\limits_{-\pi}^{\pi}|\mathcal{H}_{2}a(\theta_1, t_2)|\frac{dt_2}{2\pi}\frac{d\theta_1}{2\pi}\leq {\rm const}.
$$

So $\|\mathcal{H}a\|_{1}\leq {\rm const}$ for every
  $\mathbb{T}^2$- or  $\mathbb{T}^1$-atom $a$,  as required. $\Box$

\textbf{Proof of Proposition 2.} The first statement of the proposition  follows from Lemma 2.  To prove the second one consider a linear functional  $F\in H^1_{\mathbb{R}}(\mathbb{T}^2)^*$ such that the restriction of $F$ to the subspace generated by atoms  be zero. By Theorem 2 there is  a function $\varphi\in BMO(\mathbb{T}^2)$
such that (2) holds.
  Using  in formula (2) the atoms $q_1=(\chi+\overline{\chi})/2$ and   $q_2=(\chi-\overline{\chi})/(2i) \ (\chi\in X)$,
one deduces  that  $\widehat\varphi=0$ and therefore $F=0$.
The application of Hahn-Banach Theorem completes the proof. $\Box$

\textbf{Remark 2. }   Let $H^1_{at}$ denotes the vector space   of all function of the form $f=\sum_{j=1}^\infty\lambda_{j}a_j$ where  $a_j$ are $\mathbb{T}^2$- or $\mathbb{T}^1$-atoms, $\lambda_{j}\in \mathbb{R}$, and $\sum\limits_{j=1}^\infty|\lambda_{j}|<\infty$ endowed with the "atomic" norm
$$
\|f\|_{at}:=\inf\left\{\sum\limits_{j=1}^\infty|\lambda_{j}|:f=\sum_{j=1}^\infty\lambda_{j}a_j\right\}.
$$
By Proposition 2 there is a continuous (and dense)  embedding of $H^1_{at}$ in  $H^1_{\mathbb{R}}(\mathbb{T}^2)$. By Theorem 2 this implies a continuous embedding of $BMO(\mathbb{T}^2)$ into $(H^1_{at})^{\ast}$. In view of Hahn-Banach Theorem to prove the equality $H^1_{\mathbb{R}}(\mathbb{T}^2)=H^1_{at}$ it remains to prove  a continuous embedding of $(H^1_{at})^{\ast}$ into $BMO(\mathbb{T}^2)$.
The problem of the existence of such embedding (and therefore of an atomic decompositions for functions from  $H^1_{\mathbb{R}}(\mathbb{T}^2)$) seems to be open.

\end{document}